\documentclass[11pt,reqno]{amsart}
\usepackage[all]{xy}
\usepackage{amsmath}
\newtheorem{thm}{Theorem}%[section]
\newtheorem{lem}[thm]{Lemma}%[section]
\newtheorem{prop}[thm]{Proposition}%[section]
\theoremstyle{definition}

\theoremstyle{remark}
\theoremstyle{plain}

\newcommand{\F}{{\mathcal{F}}}
\newcommand{\OO}{{\mathcal{O}}}
\newcommand{\BB}{{\mathcal{B}}}
\newcommand{\AAA}{{\mathcal{A}}}
\newcommand{\HH}{{\mathcal{H}}}
\newcommand{\D}{{\mathcal{D}}}

\newcommand{\ZZ}{{\mathbf{Z}}}
\newcommand{\NN}{{\mathbf{N}}}
\newcommand{\CC}{{\mathbf{C}}}

\newcommand{\PP}{{\mathcal{P}}}
\newcommand{\C}{{\mathcal{C}}}

\begin{document}

%\numberwithin{equation}{section}

\title[Topological entropy]
{Topological entropy for the canonical endomorphism of Cuntz-Krieger
algebras}

\author{Florin P. Boca and Paul Goldstein}

\address{School of Mathematics, Cardiff University,
Senghennydd Road, Cardiff CF2 4YH, UK}

\address{Email of FPB: BocaFP@cardiff.ac.uk,\qquad 
Email of PG: GoldsteinP@cardiff.ac.uk}

\thanks{Research supported by an EPSRC Advanced Fellowship 
and an EPSRC Research Assistanship}
\date{June 30, 1999}

\maketitle

\setcounter{equation}{0}

Let $\Sigma$ be a finite set and let 
$A=\big( A(i,j)\big)_{i,j\in \Sigma}$ such that
$A(i,j)\in \{ 0,1\}$ and all rows and columns of $A$ are non-zero.
The Cuntz-Krieger algebra $\OO_A$ is the universal $C^*$-algebra
generated by partial isometries $S_i \neq 0$, $i\in \Sigma$,
with the property that their support projections $Q_i=S_i^* S_i$ 
and $P_i =S_iS_i^*$ satisfy the relations
$$
P_iP_j =\delta_{ij}P_i,\qquad Q_i =\sum\limits_{j\in \Sigma}
A(i,j)P_j,\qquad i,j\in \Sigma .
$$

The aim of this note is to compute the topological entropy of the 
canonical "endomorphism" $\phi_A :\OO_A \rightarrow \OO_A$, which
is the ucp (unital completely positive) map defined as
$$
\phi_A (X)=\sum\limits_{j\in \Sigma} S_j XS_j^*,\qquad 
X\in \OO_A.
$$

The map $\phi_A$ plays a crucial r\^ ole in the study of $\OO_A$ 
(\cite{CK}).
It invariates the AF-part $\F_A$ of $\OO_A$ and the abelian 
subalgebra $\D_A$ generated by $\phi_A^k (P_i)$, $i\in \Sigma$,
$k\in \NN$. The restriction $\phi_A \vert_{\D_A}$ is an isometric 
endomorphism of $\D_A$. Actually $\D_A$ identifies 
with $C(X_A)$, the commutative $C^*$-algebra of continuous functions 
on the compact space 
$$X_A =\big\{ (x_k)_{k\in \NN} \, ;\, x_k \in \Sigma ,\
A(x_k,x_{k+1})=1\big\}$$
and $\phi_A$ is the endomorphism induced on $C(X_A)$ by the 
one-sided subshift of finite type $\sigma_A$
(see \cite{CK}) defined by
$$(\sigma_A x)_k =x_{k+1},\qquad x=(x_k)_{k\in \NN} \in X_A.$$

Therefore $\phi_A$ can be regarded as a non-commutative 
generalization of the one-sided subshift of finite type
associated with the matrix $A$ and the computation
of its dynamical entropies is of some interest 
(see \cite[Page 691]{CNT}).

When $A(i,i)=1$, $i\in \Sigma$, one gets the Cuntz 
algebra $\OO_N$ where $N$ is the cardinality of $\Sigma$
(see \cite{Cu}). In this case $\phi_N =\displaystyle
\mbox{\rm \small $\sum_{j=1}^N$} S_j \cdot S_j^*$
is a genuine endomorphism (i.e. $\phi_N (XY)=\phi_N (X)
\phi_N (Y),\ X,Y\in \OO_N$) which invariates the AF-part
$\F_N =\displaystyle \mbox{\small $\bigotimes_1^{\infty}$} M_N (\CC)$
of $\OO_N$ and $\phi_N \vert_{\F_N}$ coincides with the 
noncommutative Bernoulli shift $\phi_N (X)=1\otimes X$, $X\in \F_N$. 
Furthermore, $\phi_N \vert_{\D_N}$ is the classical
one-sided Bernoulli shift.

D. Voiculescu has introduced in \cite{Vo} a notion of
topological entropy for noncommutative dynamical systems
$(A,\alpha)$, where $A$ is a unital nuclear $C^*$-algebra and
$\alpha$ an automorphism (or endomorphism) of $A$ which extends
the classical commutative topological entropy.
In the noncommutative framework partitions of unity are being 
replaced by ucp map (\cite{CNT},\cite{Vo}). As pointed out
by N. Brown (see \cite{Br}), Voiculescu's definition
carries on, with slight modifications, to the larger class of 
(not necessarily unital) exact $C^*$-algebras. 

M. Choda has computed in \cite{Cho} the topological entropy 
$\mbox{\rm ht}(\phi_N)$ of the canonical endomorphism $\phi_N$
on $\OO_N$, proving $\mbox{\rm ht}(\phi_N)=\log N$. The equality
$ht\big( \phi_A \vert _{\F_A} \big)=\log r(A)$ has 
been proved in \cite{De}. In this note 
we extend these results and compute, under a suitable definition 
for the topological entropy of a cp map, the topological entropy 
$\mbox{\rm ht}(\phi_A)$, proving

\medskip

{\bf Theorem.}
{\em If $A$ is irreducible and not a permutation matrix, then
$$
\mbox{\rm ht} (\phi_A)=\log r(A).
$$}

Here $r(A)$ denotes the spectral radius of $A$, which coincides
by Perron-Frobenius with the largest (positive) eigenvalue of $A$.

One can associate to any matrix $A=\big( A(i,j)\big)_{i,j\in \Sigma}$
with $A(i,j)\in \ZZ^+$ its dual matrix 
$A^\prime =\big( A^\prime (r,s)\big)_{r,s\in \Sigma^\prime}$
with $A^\prime (r,s)\in \{ 0,1\}$ and define
$\OO_A$ as $\OO_{A^\prime}$ (see \cite{CK}). Since $A=ST$ and
$A^\prime =TS$ for some matrices $S$ and $T$, one has
$r(A^\prime)=r(A)$. Hence the topological entropy of the canonical
endomorphism $\phi_{A^\prime}$ on $\OO_{A^\prime} =\OO_A$ equals
$r(A)$.

\bigskip

\setcounter{equation}{0}

\section{Proof of the main result}

We first recall some basic definitions from \cite{Br} and \cite{Vo}.
In the sequel $\AAA$ will be an exact $C^*$-algebra and 
$\PP \! f (\AAA)$ will denote the set of finite subsets of $\AAA$. 
For any faithful $\ast$-representation 
$\pi :\AAA \rightarrow \BB (\HH)$ one denotes by 
$CPA(\pi,\AAA)$ the set of triples $(\phi ,\psi,\BB)$, where 
$\BB$ is a finite-dimensional $C^*$-algebra and 
$\phi :\AAA \rightarrow \BB$, $\psi:\BB \rightarrow \BB (\HH)$ 
are cp maps. One also considers for any $\omega \in \PP \! f(\AAA)$ 
the completely positive $\delta$-rank
\begin{equation}
rcp (\pi,\omega;\delta)=\inf \big\{ \mbox{\rm rank} \, \BB \, ;
\, (\phi,\psi,\BB)\in CPA(\pi,\AAA),\ \| \psi \phi (a)-\pi(a) \| 
<\delta,\ a\in \omega \big\}.
\end{equation}

By an important result of E. Kirchberg and S. Wassermann exact 
$C^*$-algebras are nuclearly embeddable (see \cite{Was}). 
Hence, there exists $\pi$ faithful such that 
for all $\omega \in \PP \! f(\AAA)$ and $\delta >0$,
there is $(\phi,\psi,\BB) \in CPA(\pi,\AAA)$
with $\big\| \psi \phi (a)-\pi(a)\big\| <\delta$, $a\in \omega$. 
As noticed in \cite{Br}, $rcp (\pi ,\omega;\delta)$ is independent 
on the choice of $\pi$.

Assume also that $\Phi :\AAA \rightarrow \AAA$ is a cp map, let
$\AAA \hookrightarrow \BB (\HH)$ and define
$$\begin{array}{l}
ht (\Phi,\omega;\delta)=\limsup\limits_n n^{-1}
\log rcp \big( \omega \cup \Phi (\omega)\cup \dots
\cup \Phi^{n-1} (\omega)\, ;\, \delta \big),\\ \\
ht (\Phi ,\omega)=\sup\limits_{\delta >0}
ht (\Phi,\omega ;\delta),\qquad \quad
ht (\Phi)=\sup\limits_{\omega \in \PP \! f(\AAA)} ht (\alpha,\omega).
\end{array}$$ 

As in the case of automorphisms or endomorphisms, $ht(\Phi)$
enjoys some basic properties which are collected in the next
proposition. Proofs are similar to the corresponding ones from 
\cite{Br} and \cite{Vo}.

\medskip

\begin{prop}
$(i)$ {\em (Monotonicity)} Let $\AAA_0$ be a subalgebra of $\AAA$ 
such that $\Phi(\AAA_0)\subset \AAA_0$. Then
$$ht \big( \Phi \vert_{\AAA_0} \big) \leq ht (\Phi).$$

$(ii)$ {\em (Kolmogorov-Sinai type property)} Let $\omega_j \in \PP \!
f(\AAA)$ such that $\omega_0 \subset \omega_1 \subset \dots$ and the 
linear span of $\displaystyle \mbox{\rm \small $\bigcup_{j,k\in \NN}$}
\Phi^k (\omega_j)$ is norm dense in $\AAA$. Then
$$ht (\Phi)=\sup\limits_j ht (\Phi,\omega_j).$$

$(iii)$ {\em (Invariance to outer conjugacy)} For any
$\theta \in \mbox{\rm Aut}(\AAA)$ one has
$$ht(\theta \Phi \theta^{-1})=ht(\Phi).$$

$(iv)$ For any $k\in \NN$, $\omega \in \PP \! f (\AAA)$ and $\delta >0$
one has
$$ht (\Phi^k,\omega ;\delta)\leq k\, ht (\Phi,\omega ;\delta).$$

$(v)$ For any cp maps $\Phi_j :\AAA_j \rightarrow \AAA_j$, $j=1,2$,
one has
$$
\max \big( ht (\Phi_1),ht (\Phi_2)\big) \leq 
ht (\Phi_1 \otimes \Phi_2) \leq ht (\Phi_1)+ht (\Phi_2) .
$$
\end{prop}

\smallskip

Next we turn to the Cuntz-Krieger $C^*$-algebra $\OO_A$ associated with 
an irreducible, non-permutation matrix $A$ with entries 
in $\{ 0,1\}$. For any $k$-tuple 
$\mu =(\mu_1,\dots,\mu_k)$, $\mu_j \in \Sigma$,
we denote $o(\mu)=\mu_1$, $t(\mu)=\mu_k$, $S_\mu =S_{\mu_1} \dots 
S_{\mu_k}$, $S_e =I$, $o(e)=t(e)=I$ ($e$ denotes the empty word)
and $Q_\mu =S_\mu^* S_\mu$. For $\mu =(\mu_1,\dots,\mu_k),
\nu =(\nu_1,\dots,\nu_l)$, $\mu_i,\nu_j \in \Sigma$ we denote
$\mu \nu =(\mu_1,\dots,\mu_k,\nu_1,\dots,\nu_l)$.
The number of elements of a finite set $F$ is denoted by $\# F$.
We set $A(\mu)=1$ for $k=1$ and 
$$A(\mu)=A(\mu_1,\mu_2)A(\mu_2,\mu_3)\dots A(\mu_{k-1} ,\mu_k)
\qquad \mbox{\rm for $k\geq 2$}.
$$

Then, for $\mu$, $\nu$ with $\vert \mu \vert =\vert \nu \vert$ one has
$$\begin{array}{l}
S_\mu^* S_\nu =\delta_{\mu \nu} Q_\mu =\delta_{\mu \nu} A(\mu)
Q_{t(\mu)},\\ \\
Q_\eta S_\alpha =A\big( \eta o(\alpha)\big) S_\alpha
\qquad \mbox{\rm for $\vert \alpha \vert \geq 1$}.\end{array}$$

In particular $S_\mu \neq 0$, $\vert \mu \vert \geq 1$, is 
equivalent to $A(\mu) \neq 0$. It is clear that the number of
elements of 
$$
L(k)=\big\{ \mu \, ;\, \vert \mu \vert =k,\ S_\mu \neq 0\big\}
$$
equals
\begin{equation}
\begin{array}{rl} w(k) & =\# \big\{ (i_1,\dots,i_k)\, ;\, 
i_j \in \Sigma ,\ A(i_1,i_2)A(i_2,i_3) \dots A(i_{k-1},i_k)=1\big\}
\\   \\ \qquad & =\displaystyle 
\sum\limits_{i_1,\dots ,i_k\in \Sigma}
A(i_1,i_2)A(i_2,i_3) \dots A(i_{k-1},i_k) =
\sum\limits_{i,j\in \Sigma} A^{k-1} (i,j) \\   \displaystyle
\qquad & =\langle A^{k-1} e,e\rangle, \qquad \qquad 
\mbox{\rm where}\ e=\begin{pmatrix} 1 \\ \vdots \\ 1 \end{pmatrix} .
\end{array}
\end{equation}

Note also that if $r(A)$ denotes the spectral radius of $A$, then
$$\| A^{k-1} \| \leq \sum\limits_{i,j\in \Sigma} A^{k-1} (i,j)
=w(k)\leq \| A^{k-1} \| \cdot \| e\|_2^2 =\# \Sigma \cdot 
\| A^{k-1} \| ,$$
which provides 
\begin{equation}
\lim\limits_k k^{-1} \log w(k)=\lim\limits_k k^{-1} \log \| A^{k-1} \|
=\log r(A).
\end{equation}

We consider now a certain embedding of the Cuntz-Krieger algebra 
$\OO_A$ into $M_{w(m)} (\CC)\otimes \OO_A$. For each $m\geq 1$ we index
the canonical matrix unit of $M_{w(m)} (\CC)$ as
$\{ e_{\mu \nu} \}_{\mu,\nu \in L(m)}$ and define a map 
$\rho_m :\OO_A \rightarrow M_{w(m)} (\CC)\otimes \OO_A$ by
\begin{equation}
\rho_m (X)=\sum\limits_{\mu,\nu \in L(m)} e_{\mu \nu} \otimes
S_\mu^* XS_\nu.
\end{equation}

Since $\displaystyle \mbox{\small $\sum_{\mu \in L(m)}$} 
S_\mu S_\mu^* =\mbox{\small $\sum_{\vert \mu \vert =m}$} 
S_\mu S_\mu^*=I$, it is easily seen that $\rho_m$ is a 
$\ast$-morphism. Moreover,
since $\OO_A$ is simple, it follows that $\rho_m :\OO_A \rightarrow 
\rho_m (\OO_A) \subset M_{w(m)} (\CC)\otimes \OO_A$ is a 
$\ast$-isomorphism. The map $\rho_m$ is not unital in general. We only 
have $\rho_m(1)=\displaystyle \mbox{\small $\sum_{\mu \in L(m)}$}
e_{\mu \mu} \otimes Q_\mu =\mbox{\small $\sum_{\mu \in L(m)}$}
e_{\mu \mu} \otimes Q_{t(\mu)}$. 

However, for Cuntz algebras this map is unital, multiplicative 
and onto, providing an explicit isomorphism between $\OO_N$ and
$M_{N^m} (\CC) \otimes \OO_N$.  To see that $\rho_m$ is onto note that
for $\vert \mu_0 \vert=\vert \nu_0 \vert =m$ we have
$\rho_m (S_{\mu_0} S_{\nu_0}^*) =e_{\mu_0 \nu_0} \otimes I$ and
$\rho_m (S_{\mu_0} S_{\nu_0} S_{\mu_0}^*)=e_{\mu_0 \mu_0} \otimes
S_{\nu_0}$. For $m=1$ and $N=2$ this map was used by M. D. Choi 
(see \cite{Ch}) to prove the isomorphism between $M_2(\CC)\otimes
\OO_2$ and $\OO_2$.

We return to the general case and note that for all $l\geq 1$
\begin{equation}
\phi_A^l (X)=\sum\limits_{\vert \eta \vert =l} S_\eta XS_\eta^* =
\sum\limits_{\eta \in L(l)} S_\eta XS_\eta^*,\qquad X\in \OO_A.
\end{equation}

\medskip

\begin{lem}
Let $n\geq 1$ and assume that 
$\vert \beta \vert \leq \vert \alpha \vert \leq n_0$ and
$m\geq n+n_0$. Then, for all $i\in \Sigma$ and all
$l\in \{ 0,1,\dots,n-1\}$ one has
$$\rho_m \phi_A^l (S_\alpha P_i S_\beta^*) =
\begin{cases} \displaystyle
\sum\limits_{\vert \mu \vert =\vert \alpha \vert 
-\vert \beta \vert} X(\mu) \otimes S_\mu & \mbox{\rm if
$\vert \beta \vert <\vert \alpha \vert$},\\ \displaystyle
\sum\limits_{j\in \Sigma} X_j \otimes Q_j & \mbox{\rm if $\vert
\beta \vert =\vert \alpha \vert$},\end{cases}$$
for some partial isometries $X(\mu)=X\big( \vert \alpha \vert ,\vert 
\beta \vert ,i,l,\mu \big)$ and respectively $X_j =X\big( \vert
\alpha \vert ,i,l,j\big)$.
\end{lem}

{\sl Proof.}
From (5) and (4) we get
\begin{equation}
\begin{array}{rl} 
\rho_m \phi_A^l (S_\alpha P_i S_\beta^*) & \displaystyle =
\sum\limits_{\eta \in L(l)} \rho_m (S_{\eta \alpha} P_i S_{\eta \beta}^*)
=\sum\limits_{\eta \in L(l)} \
\sum\limits_{\mu,\nu \in L(m)}
e_{\mu \nu} \otimes S_\mu^* S_\eta S_\alpha P_i S_\beta^* S_\eta^* S_\nu
\\ (with \
\mu =\eta \mu^\prime,\ \nu =\eta \nu^\prime) & \displaystyle  
=\sum\limits_{\vert \eta \vert =l} \
\sum_{\substack{\vert \mu^\prime \vert =\vert \nu^\prime \vert
=m-l\\ \eta \mu^\prime ,\eta \nu^\prime \in L(m)}}
e_{\eta \mu^\prime,\eta \nu^\prime} \otimes 
S_{\mu^\prime}^* Q_\eta S_\alpha P_i S_\beta^* Q_\eta S_{\nu^\prime}.
\end{array}
\end{equation}

For $\vert \beta \vert =\vert \alpha \vert \geq 1$ this yields
$$\begin{array}{rl} 
\rho_m \phi_A^l (S_\alpha P_i S_\beta^*)
& \displaystyle =
\sum\limits_{\vert \eta \vert =l} \
\sum_{\substack{\vert \mu^\prime \vert =\vert \nu^\prime \vert =m-l\\ 
\eta \mu^\prime,\eta \nu^\prime \in L(m)}}
A\big( \eta o(\alpha)\big) A\big( \eta o(\beta)\big) 
e_{\eta \mu^\prime,\eta \nu^\prime} \otimes
S_{\mu^\prime}^* S_\alpha P_i S_\beta^* S_{\nu^\prime} \\
(with \ \mu^\prime=\alpha \mu^{\prime \prime},
\ \nu^\prime =\beta \nu^{\prime \prime}) & \displaystyle
=\sum\limits_{\vert \eta \vert =l} \
\sum_{\substack{\vert \mu^{\prime \prime} \vert =\vert
\nu^{\prime \prime} \vert =m-l-\vert \alpha \vert\\ 
\eta \alpha \mu^{\prime \prime},\eta \beta \nu^{\prime \prime}
\in L(m)}} A\big( \eta o(\alpha)\big) A\big( \eta o(\beta)\big) 
e_{\eta \alpha \mu^{\prime \prime},\eta \beta \nu^{\prime \prime}}
\otimes S_{\mu^{\prime \prime}}^* Q_\alpha P_i 
Q_\beta S_{\nu^{\prime \prime}} \\ & \displaystyle
=\sum\limits_{\vert \eta \vert =l} \
\sum_{\substack{\vert \mu^{\prime \prime} \vert =
\vert \nu^{\prime \prime} \vert =m-l-\vert \alpha \vert\\ 
\eta \alpha \mu^{\prime \prime},\eta \beta \nu^{\prime \prime}
\in L(m)}} A\big( \eta \alpha o(\mu^{\prime \prime})\big)
A\big( \eta \beta o(\nu^{\prime \prime})\big)  
e_{\eta \alpha \mu^{\prime \prime},\eta \beta \nu^{\prime \prime}}
\otimes S_{\mu^{\prime \prime}}^* P_i S_{\nu^{\prime \prime}}
\\  & \displaystyle =\sum\limits_{\vert \eta \vert =l} \
\sum_{\substack{\vert \mu^{\prime \prime} \vert =
\vert \nu^{\prime \prime} \vert =m-l-\vert \alpha \vert \\
o(\mu^{\prime \prime})=o(\nu^{\prime \prime})=i \\ 
\eta \alpha \mu^{\prime \prime},\eta \beta \nu^{\prime \prime}
\in L(m)}} A(\eta \alpha i)A(\eta \beta i) 
e_{\eta \alpha \mu^{\prime \prime},\eta \beta \nu^{\prime \prime}}
\otimes S_{\mu^{\prime \prime}}^* S_{\nu^{\prime \prime}} \\
(with \ \nu^{\prime \prime}=\mu^{\prime \prime}) & \displaystyle
=\sum\limits_{\vert \eta \vert =l} \
\sum_{\substack{\vert \mu^{\prime \prime} \vert =
m-l-\vert \alpha \vert,\, o(\mu^{\prime \prime}) =i\\ 
\eta \alpha \mu^{\prime \prime},\eta \beta \mu^{\prime \prime}
\in L(m)}} A(\eta \alpha i)A(\eta \beta i) A(\mu^{\prime \prime})
e_{\eta \alpha \mu^{\prime \prime},\eta \beta \mu^{\prime \prime}}
\otimes Q_{t(\mu^{\prime \prime})} \\  & \displaystyle 
=\sum\limits_{\vert \eta \vert =l} \
\sum_{\substack{\vert \mu^{\prime \prime} \vert =
m-l-\vert \alpha \vert,\, o(\mu^{\prime \prime}) =i\\ 
\eta \alpha \mu^{\prime \prime},\eta \beta \mu^{\prime \prime}
\in L(m)}} 
e_{\eta \alpha \mu^{\prime \prime},\eta \beta \mu^{\prime \prime}}
\otimes Q_{t(\mu^{\prime \prime})} \\  & \displaystyle 
=\sum\limits_{j\in \Sigma} X_j \otimes Q_j,
\end{array}$$
where
$$X_j =\sum\limits_{\vert \eta \vert =l} \
\sum_{\substack{\vert \mu^{\prime \prime} \vert =
m-l-\vert \alpha \vert\\ o(\mu^{\prime \prime})=i,\,
t(\mu^{\prime \prime})=j\\
\eta \alpha \mu^{\prime \prime},\eta \beta \mu^{\prime \prime} \in L(m)}}
e_{\eta \alpha \mu^{\prime \prime},\eta \beta \mu^{\prime \prime}}$$
are partial isometries for all $j\in \Sigma$.

For $\beta =\alpha =e$ a similar computation yields
$\rho_m \phi_A^l (P_i)=\displaystyle 
\mbox{\small $\sum_{j\in \Sigma}$} X_j \otimes Q_j$, with
$$X_j =\sum\limits_{\vert \eta \vert =l} \
\sum_{\substack{\vert \mu^\prime \vert =m-l,\, 
o(\mu^\prime)=i,\, t(\mu^\prime)=j\\
\eta \mu^\prime \in L(m)}} e_{\eta \mu^\prime,\eta \mu^\prime}.$$

For $1\leq \vert \beta \vert <\vert \alpha \vert$ equality (6) yields
$$\begin{array}{rl}
\rho_m \phi_A^l (S_\alpha P_i S_\beta^*) & \displaystyle
=\sum\limits_{\vert \eta \vert =l} \
\sum_{\substack{\vert \mu^\prime \vert =\vert \nu^\prime \vert
=m-l \\ \eta \mu^\prime,\eta \nu^\prime \in L(m)}}
A\big( \eta o(\alpha)\big) A\big( \eta o(\beta)\big)
e_{\eta \mu^\prime, \eta \nu^\prime} \otimes
S_{\mu^\prime}^* S_\alpha P_i S_\beta^* S_{\nu^\prime} \\
(with \ \mu^\prime =\alpha \mu^{\prime \prime},
\ \nu^\prime =\beta \nu^{\prime \prime}) & \displaystyle =
\sum\limits_{\vert \eta \vert =l} \
\sum_{\substack{\vert \mu^{\prime \prime} \vert =
m-l-\vert \alpha \vert\\
\vert \nu^{\prime \prime} \vert =m-l-\vert \beta \vert \\
\eta \alpha \mu^{\prime \prime},\eta \beta \nu^{\prime \prime}
\in L(m)}} A\big( \eta o(\alpha)\big) A\big( \eta o(\beta)\big)
e_{\eta \alpha \mu^{\prime \prime},\eta \beta \nu^{\prime \prime}}
\otimes S_{\mu^{\prime \prime}}^* Q_\alpha P_i 
Q_\beta S_{\nu^{\prime \prime}} \\ & \displaystyle
=\sum\limits_{\vert \eta \vert =l} \
\sum_{\substack{\vert \mu^{\prime \prime} \vert =
m-l-\vert \alpha \vert \\
\vert \nu^{\prime \prime} \vert =m-l-\vert \beta \vert \\
\eta \alpha \mu^{\prime \prime},\eta \beta \nu^{\prime \prime}
\in L(m)}} A\big( \eta \alpha o(\mu^{\prime \prime})\big)
A\big( \eta \beta o(\nu^{\prime \prime})\big) 
e_{\eta \alpha \mu^{\prime \prime}, \eta \beta \nu^{\prime \prime}}
\otimes S_{\mu^{\prime \prime}}^* P_i S_{\nu^{\prime \prime}} \\
& \displaystyle =\sum\limits_{\vert \eta \vert =l} \
\sum_{\substack{\vert \mu^{\prime \prime} \vert =m-l-\vert \alpha \vert,
\, o(\mu^{\prime \prime})=i \\
\vert \nu^{\prime \prime} \vert = m-l-\vert \beta \vert,\,
o(\nu^{\prime \prime})=i\\
\eta \alpha \mu^{\prime \prime},\eta \beta \nu^{\prime \prime}
\in L(m)}} A(\eta \alpha i)A(\eta \beta i)
e_{\eta \alpha \mu^{\prime \prime},\eta \beta \nu^{\prime \prime}}
\otimes S_{\mu^{\prime \prime}}^* S_{\nu^{\prime \prime}} \\
(with \ 
\nu^{\prime \prime}=\mu^{\prime \prime} \mu) & \displaystyle =
\sum\limits_{\vert \eta \vert =l} \
\sum_{\substack{\vert \mu^{\prime \prime} \vert =m-l-\vert \alpha \vert,
\, o(\mu^{\prime \prime})=i\\ 
\vert \mu \vert =\vert \alpha \vert -\vert \beta \vert \\
\eta \alpha \mu^{\prime \prime},
\eta \beta \mu^{\prime \prime} \! \mu \in L(m)}}
A\big( \eta \alpha \mu^{\prime \prime} o(\mu)\big)
A(\eta \beta i) e_{\eta \alpha \mu^{\prime \prime},\eta \beta
\mu^{\prime \prime} \mu} \otimes S_\mu \\ & \displaystyle =
\sum\limits_{\vert \eta \vert =l} \
\sum_{\substack{\vert \mu^{\prime \prime} \vert=m-l-\vert \alpha \vert,
\, o(\mu^{\prime \prime})=i\\ 
\vert \mu \vert =\vert \alpha \vert -\vert \beta \vert \\
\eta \alpha \mu^{\prime \prime},
\eta \beta \mu^{\prime \prime} \! \mu \in L(m)}}
e_{\eta \alpha \mu^{\prime \prime},\eta \beta
\mu^{\prime \prime} \mu} \otimes S_\mu \\  & \displaystyle =
\sum\limits_{\vert \mu \vert =\vert \alpha \vert -\vert \beta \vert} 
X(\mu)\otimes S_\mu,
\end{array}$$
where
$$X(\mu)=\sum\limits_{\vert \eta \vert =l} \
\sum_{\substack{\vert \mu^{\prime \prime} \vert=m-l-\vert \alpha \vert,
\, o(\mu^{\prime \prime})=i\\
\eta \alpha \mu^{\prime \prime},\eta \beta \mu^{\prime \prime} \! \mu
\in L(m)}} e_{\eta \alpha \mu^{\prime \prime},
\eta \beta \mu^{\prime \prime} \mu}$$
are partial isometries for all $\mu \in L(\vert \alpha \vert -
\vert \beta \vert)$. One plainly checks that for $\beta =e$,
$\vert \alpha \vert \geq 1$, the formula
$\rho_m \phi_A^l (S_\alpha P_i)=
\displaystyle \mbox{\small $\sum_{\vert \mu \vert =\vert \alpha \vert}$}
X(\mu)\otimes S_\mu$ holds for the $X(\mu)$ above which corresponds
to $\beta =e$.$\qquad \qquad \qed$
\medskip

For any $k\geq 1$ we put
$$\omega (k)=\big\{ S_\alpha P_i S_\beta^* \, ;\, \vert \beta \vert 
\leq \vert \alpha \vert \leq k,\ i\in \Sigma \big\}.$$

\medskip

\begin{prop}
For all $n_0 \geq 1$ and $\delta >0$ one has
$$\limsup\limits_n n^{-1} \log rcp \left( \omega(n_0)
\cup \phi_A \big( \omega(n_0)\big) \cup \dots \cup \phi_A^{n-1} 
\big( \omega(n_0)\big)\, ; \, \delta \right) \leq \log r(A).$$
\end{prop}

{\sl Proof.}
For $n\geq 1$ we let $m=m(n)=n+n_0$. Since $\OO_A$ is nuclear, 
there exists $\big( \phi_0,\psi_0,M_{m_0} (\CC)\big) \in 
CPA(id_A,\OO_A)$, that is
$$
\xymatrix{ \mbox{\footnotesize $\OO_A$} \ar@{->} [rr]^{id_{\OO_A}} 
\ar@{->} [dr] _{\phi_0} & & \mbox{\footnotesize $\OO_A$} \\
& \mbox{\footnotesize $M_{m_0} (\CC)$} \ar@{->} [ur] _{\psi_0} & }$$
such that
\begin{equation}
\big\| \psi_0 \phi_0 (Q_j)-Q_j\big\| +
\big\| \psi_0 \phi_0 (S_\gamma)-S_\gamma \big\|
<\frac{\delta}{\max \big( \# \Sigma,w(n_0)\big)} \qquad 
\mbox{\rm for all $\gamma \in L(n_0)$ and $j\in \Sigma$.}
\end{equation}

Consider $\BB=M_{w(m)} (\CC )\otimes M_{m_0} (\CC)$ and let $\HH$ be
a Hilbert space on which $\OO_A$ acts faithfully. The $\ast$-isomorphism
$\rho_m^{-1}:\rho_m(\OO_A) \rightarrow \OO_A$ extends to a cp map
$\Psi_m :M_{w(m)} (\CC)\otimes \OO_A \rightarrow \BB (\HH)$
with $\| \Psi_m \| =1$. We consider the cp maps
$\phi =(id \otimes \phi_0)\rho_m :\OO_A \rightarrow \BB$ and 
$\psi =\Psi_m (id \otimes \psi_0):\BB \rightarrow \BB (\HH)$; 
see the following diagram
$$\xymatrix@M=4pt{
\mbox{\footnotesize $\OO_A$} \ar@{->} [rr] ^(0.45){\rho_m}  
\ar@{->} [ddrrr] _(0.4){\phi} & & \ar @{} [d] |{=} 
\mbox{\footnotesize $\rho_m(\OO_A)$}
\ar@{->} [ddr] ^{id \otimes \phi_0} \ar@{->} [rr] ^{id_{\rho_m(\OO_A)}}
 &  & 
\ar @{} [dr] < 8pt> |{=} \mbox{\footnotesize $\rho_m(\OO_A)$}
 \ar@{->} [rr] ^{\rho_m^{-1}} 
\ar@{^{(}->} [d] &  & \mbox{\footnotesize $\BB (\HH)$}  \\ 
&   &  &   & \mbox{\footnotesize $M_{w(m)} (\CC) \otimes \OO_A$}
\ar @{} [d]  |{=}  \ar@{-->} [urr] _{\Psi_m} &  & \\
         &     &   & 
\mbox{\footnotesize $\BB=M_{w(m)} (\CC)\otimes M_{m_0} (\CC)$}  
\ar@{->} [ur] _(0.55){id \otimes \psi_0} 
\ar@/_3pc/ @{->} [uurrr] _(0.65){\psi} &   &   &  }$$

Let $a=S_\alpha P_i S_\beta^* \in \omega (n_0)$. By the previous lemma
there exist partial isometries $X(\mu)=X(a,l,\mu )$ if $\vert \beta 
\vert <\vert \alpha \vert$ and $X_j =X(a,l,j)$ if
$\vert \alpha \vert =\vert \beta \vert$ such that
\begin{equation}
\rho_m \phi_A^l (a) =\begin{cases} \displaystyle
\sum\limits_{\vert \mu\vert =\vert \alpha \vert -
\vert \beta \vert } X(\mu)\otimes S_\mu &
\mbox{\rm for $\vert \beta \vert <\vert \alpha \vert$,} \\
\displaystyle \sum\limits_{j\in \Sigma} X_j \otimes Q_j
& \mbox{\rm for $\vert \beta \vert =\vert \alpha \vert$.}
\end{cases} \end{equation}

From (8) and (7) we gather
$$\begin{array}{l} 
\big\| \psi \phi \big( \phi_A^l (a)\big) -\phi_A^l (a)\big\| 
=\big\| \Psi_m (id \otimes \psi_0 \phi_0 )\big( \rho_m \phi_A^l (a)\big)
-\phi_A^l(a)\big\| \\ \\  \displaystyle \qquad 
=\big\| \Psi_m (id \otimes \psi_0 \phi_0)
\big( \rho_m \phi_A^l (a)\big) -\Psi_m \big( \rho_m \phi_A^l (a)\big) 
\big\| \leq \big\| (id \otimes \psi_0 \phi_0)\big( 
\rho_m \phi_A^l (a)\big) -\rho_m \phi_A^l (a)\big\| \\ \\ 
\qquad =\begin{cases} \displaystyle 
\Bigl\| \sum\limits_{\vert \mu \vert =\vert
\alpha \vert -\vert \beta \vert } X(\mu) \otimes 
\big( \psi_0 \phi_0 (S_\mu)-S_\mu \big) \Bigr\| & \mbox{\rm for 
$\vert \beta \vert <\vert \alpha \vert$,} \\  \displaystyle
\Bigl\| \sum\limits_{j\in \Sigma} X_j \otimes \big( 
\psi_0 \phi_0 (Q_j)-Q_j\big) \Bigr\| & \mbox{\rm for
$\vert \beta \vert =\vert \alpha \vert$,} \end{cases} \\  \\  
\displaystyle \qquad 
< \max \big( \# \Sigma,w(n_0)\big) \cdot 
\frac{\delta}{\max \big( \# \Sigma, w(n_0)\big)}
\, =\delta.
\end{array}$$

Therefore
$$rcp \left( \omega(n_0)\cup \phi_A \big( \omega(n_0)\big) \cup
\dots \cup \phi_A^{n-1}\big( \omega(n_0)\big) \, ;\, \delta \right) 
\leq m_0 w(m) =m_0 w(n+n_0),$$
which we combine with (3) to get
$$\limsup\limits_n n^{-1} \log \left( \omega(n_0)
\cup \dots \cup \phi_A^{n-1} \big( \omega(n_0)\big)  
\, ;\, \delta\right) \leq \limsup\limits_n
n^{-1} \log w(n) =\log r(A).\ \qed$$

\medskip

{\sl Proof of the main result.}
Since $\omega_n =\omega(n)\cup \omega(n)^*$ is an increasing sequence
of finite subsets of $\OO_A$ and $\mbox{\rm span} \, \displaystyle
\mbox{\rm \small $\bigcup_n$} \omega_n$ is dense in the
uniform norm in $\OO_A$, Proposition 1$\, (ii)$ and Proposition 3
provide
\begin{equation}
ht(\phi_A)\leq \log r(A).
\end{equation}

For the opposite inequality, denote 
$\theta_A=\phi_A \vert_{\D_A =C(X_A)}$. By Proposition 1$\, (i)$ 
$ht(\phi_A)\geq ht(\theta_A)$. Within the framework of \cite{CNT},
let $\sigma$ be a probability measure on $X_A$ such that
$\sigma \theta_A=\sigma$. For any finite-dimensional algebra $M$
and any ucp map $\gamma :M\rightarrow C(X_A)$, Proposition III.6
in \cite{CNT} provides $H_\sigma (\gamma,\theta_A \gamma,\dots,
\theta_A^{n+m-1} \gamma) \leq H_\sigma (\gamma,\theta_A \gamma,
\dots,\theta_A^{n-1} \gamma)+H_\sigma (\gamma,\theta_A \gamma,
\dots,\theta_A^{n-1} \gamma)$ for all $m,n\geq 1$, hence
$$h_{\sigma,\theta_A} (\gamma) =\lim\limits_n n^{-1}
H_\sigma (\gamma,\theta_A \gamma,\dots ,\theta_A^{n-1} \gamma)$$
exists. Let $h_\sigma (\theta_A)$ be the supremum of
$h_{\sigma,\theta_A}(\gamma)$ over all such $M$ and $\gamma$.
Arguing as in \cite[Prop.4.8]{Vo} one proves that for any $\gamma
:M\rightarrow C(X_A)$ as above and any $\varepsilon >0$, there
exist $\omega \in \PP \! f\big( C(X_A)\big)$ and $\delta >0$ such that
$$H_\sigma (\gamma,\theta_A \gamma ,\dots,\theta_A^{n-1} \gamma) \leq
n\varepsilon +\log rcp \left( \omega \cup \theta_A(\omega)\cup 
\dots \cup \theta_A^{n-1} (\omega);\delta \right),$$
hence
\begin{equation}
h_{\sigma,\theta_A} (\gamma) \leq ht(\theta_A).
\end{equation}

If $\PP$ is a finite partition into time-zero cylinder sets,
$\PP_n =\PP \vee \sigma_A^{-1} \PP \vee \dots \vee \sigma_A^{-(n-1)}
\PP$, $\C$ is the abelian finite-dimensional $C^*$-algebra
generated by $\PP_n$ and $\gamma=i_{\C}$ the natural inclusion of 
$\C$ into $C(X_A)$, then
\begin{equation}
-\sum\limits_{E\in \PP_n} \sigma (\chi_E) \log \sigma(\chi_E) =
S\big( \sigma \vert_{\PP_n} \big) =
H_\sigma (\gamma,\theta_A \gamma ,\dots,\theta_A^{n-1} \gamma),
\end{equation}
the last equality following from \cite[Remark III.5.2]{CNT}.
From (10) and (11) it follows that the classical measurable
entropy $H_\sigma (\sigma_A)$ is $\leq ht(\theta_A)$. 
In the case when $\sigma$ is the probability measure defined by
a probability eigenvector of $A$ it is well-known
(see e.g. \cite{Pet}) that $H_\sigma (\sigma_A)=\log r(A)$.
Hence one has $ht(\phi_A)\geq ht(\theta_A)\geq \log r(A)$,
which completes the proof.$\qquad \qquad \qed$

\bigskip

\bigskip

\bibliographystyle{amsplain}

\end{document}